\documentclass[12pt,leqno]{article}
\usepackage{amsfonts}
\usepackage{amsmath, amsthm, amsfonts, amssymb, color}
\usepackage{mathrsfs}
\theoremstyle{definition}

\newcommand{\scr}[1]{\mathscr #1}
\definecolor{wco}{rgb}{0.5,0.2,0.3}

\numberwithin{equation}{section} \theoremstyle{remark}

\newcommand{\ua}{\uparrow}

\title{{\bf Generalized Curvature Condition for Subelliptic Diffusion Processes}\footnote{Supported in
 part by  Lab. Math. Com. Sys., NNSFC(11131003), SRFDP and the Fundamental Research Funds for the Central Universities.}
}
\author{
{\bf   Feng-Yu Wang}\\
\footnotesize{School of Mathematical Sciences,
Beijing Normal
University, Beijing 100875, China}\\
 \footnotesize{Department of Mathematics,
Swansea University, Singleton Park, SA2 8PP, United Kingdom}\\ \footnotesize{wangfy@bnu.edu.cn, F.-Y.Wang@swansea.ac.uk}}
\begin{document}
\def\R{\mathbb R}  \def\ff{\frac} \def\ss{\sqrt} \def\B{\mathbf
B}
\def\N{\mathbb N} \def\kk{\kappa} \def\m{{\bf m}}
\def\dd{\delta} \def\DD{\Delta} \def\vv{\varepsilon} \def\rr{\rho}
\def\<{\langle} \def\>{\rangle} \def\GG{\Gamma} \def\gg{\gamma}
  \def\nn{\nabla} \def\pp{\partial} \def\EE{\scr E}
\def\d{\text{\rm{d}}} \def\bb{\beta} \def\aa{\alpha} \def\D{\scr D}
  \def\si{\sigma} \def\ess{\text{\rm{ess}}}
\def\beg{\begin} \def\beq{\begin{equation}}  \def\F{\scr F}
\def\Ric{\text{\rm{Ric}}} \def\Hess{\text{\rm{Hess}}}
\def\e{\text{\rm{e}}} \def\ua{\underline a} \def\OO{\Omega}  \def\oo{\omega}
 \def\tt{\tilde} \def\Ric{\text{\rm{Ric}}}
\def\cut{\text{\rm{cut}}} \def\P{\mathbb P} \def\ifn{I_n(f^{\bigotimes n})}
\def\C{\scr C}      \def\aaa{\mathbf{r}}     \def\r{r}
\def\gap{\text{\rm{gap}}} \def\prr{\pi_{{\bf m},\varrho}}  \def\r{\mathbf r}
\def\Z{\mathbb Z} \def\vrr{\varrho} \def\ll{\lambda}
\def\L{\scr L}\def\Tt{\tt} \def\TT{\tt}\def\II{\mathbb I}
\def\i{{\rm in}}\def\Sect{{\rm Sect}}\def\E{\mathbb E} \def\H{\mathbb H}
\def\M{\scr M}\def\Q{\mathbb Q} \def\texto{\text{o}} \def\LL{\Lambda}
\def\Rank{{\rm Rank}} \def\B{\scr B} \def\i{{\rm i}} \def\HR{\hat{\R}^d}
\def\to{\rightarrow}\def\l{\ell}
\def\8{\infty}

\maketitle

\begin{abstract}  By using a general version of curvature condition, derivative inequalities are established for a large class of subelliptic diffusion semigroups. As applications,  the Harnack/cost-entropy/cost-variance inequalities for the diffusion semigroups, and the Poincar\'e/log-Sobolev inequalities for the associated Dirichlet forms in the symmetric case,  are derived. Our results largely generalize  and partly improve the corresponding ones   obtained recently in \cite{BB}.  \end{abstract} \noindent
 \vskip 2cm

\section{Introduction}

It is well known that Bakry-Emery's curvature and curvature-dimension conditions have played  crucial roles in the study of elliptic diffusion processes. When
the diffusion operator is merely subelliptic, this condition is however no longer available.
Recently, in order to study subelliptic diffusion processes, a nice generalized curvature-dimension condition was introduced and applied in \cite{BB,BBG,BG1}, so that many important results derived in the elliptic setting have been extended to subelliptic diffusion processes with generators of type
$$L:=\sum_{i=1}^n X_i^2+X_0$$ for smooth vector fields $\{X_i: 0\le n\le 1\}$ on a differential manifold such that $\{X_i, \nn_{X_i}X_j: 1\le i\le n \}$ spans the tangent space (see comment (a) below). See also \cite{BBBC, BGM, DM, Li} and the references within for the study of the heat semigroup generated by the Kohn-Laplacian on Heisenberg groups by other means.  Stimulated by \cite{BB}, in this paper we aim to introduce a new generalized curvature condition to study more general subelliptic diffusion processes.

Let $M$ be a connected differential manifold, and let $L$ be given above
for some $C^2$-vector fields  $\{X_i\}_{i=1}^d$ and a $C^1$-vector field  $X_0$. The square field (or carr\'e du champ) for $L$ is a symmetric bilinear differential form given by
$$\GG(f,g)=\sum_{i=1}^n (X_if)(X_ig),\ \ f,g\in C^1(M).$$ Obviously, $\GG$ satisfies
$$\GG(f):=\GG(f,f)\ge 0,\ \GG(fg,h)= g\GG(f,h)+f\GG(g,h),\ \GG(\phi\circ f, g)= (\phi'\circ f)\GG(f,g)$$ for any $f,g,h\in C^1(M)$ and  $\phi\in C^1(\R).$
From now on,
a symmetric bilinear differential form $\bar\GG$ satisfying these properties is called a diffusion square field. If moreover  for any
$x\in M$ and $f\in C^1(M)$, $\bar\GG(f)(x)=0$ implies $(\d f)(x)=0$, we call $\bar\GG$ elliptic or non-degenerate.

For any $C^2$-diffusion square field $\bar \GG$ (i.e. $\bar\GG(f,g)\in C^2(M)$ for $f,g\in C^\infty(M)$), we define the associated Bakry-Emery curvature operator w.r.t. $L$ by
$$\bar\GG_2(f)= \ff 1 2L\bar \GG(f)-\bar\GG(f,Lf),\ \ f\in C^3(M).$$
Then the generalized curvature-dimension condition introduced in \cite{BG1} reads
\beq\label{GCD} \GG_2(f)+r\GG^Z_2(f) \ge \ff{(Lf)^2} d +\Big(\rr_1-\ff\kk r\Big)\GG(f) +\rr_2\GG^Z(f),\ \ f\in C^2(M), r>0,\end{equation}
where $\rr_2>0,\kk\ge 0,\rr_1\in\R$ and $d\in (0,\infty]$ are constants, and
$\GG^Z$ is a
$C^2$-diffusion square field such that $\GG+\GG^Z$ is elliptic and
\beq\label{EXX}\GG(\GG^Z(f), f)= \GG^Z(\GG(f),f),\ \ f\in C^\infty(M)\end{equation} holds.  When $\GG^Z=0$, (\ref{GCD}) reduces back to the Bakry-Emery curvature-dimension condition \cite{BE}, and when $d=\infty$ it becomes the following generalized curvature condition
\beq\label{GCD'} \GG_2(f)+r\GG^Z(f) \ge    \Big(\rr_1-\ff\kk r\Big)\GG(f) +\rr_2\GG^Z(f),\ \ f\in C^2(M), r>0.\end{equation}
Using (\ref{EXX}) and (\ref{GCD'}) for symmetric subelliptic operators, Poincar\'e inequality for the associated Dirichlet form, and
the Harnack inequality and the  log-Sobolev inequality
(for, however,  an enlarged Dirichlet form given by $\GG+\GG^Z$) for the associated diffusion semigroup, and the HWI inequality (where
the energy part is  given by the enlarged Dirichlet form) are investigated in \cite{BB}.

The main purpose of this paper is to introduce a  general version of the curvature condition to derive better inequalities for  more general subelliptic diffusion semigroups. The necessity of our study is based on the following three observations:
\beg{enumerate} \item[(a)] Condition (\ref{GCD'}) is not available if the family $\{X_i, \nn_{X_i}X_j\}_{1\le i,j\le n}$ does not span the tangent space $TM$, where $\nn$ is the Levi-Civita connection w.r.t. a Riemannian metric. Indeed, (\ref{GCD'}) implies  that the diffusion square field
$$\tt\GG(f):= \GG(f)+ \sum_{i,j=1}^n\big((\nn_{X_i} X_j) f\big)^2,\  \ f\in C^1(M)$$ is elliptic. To see this, let $x\in M$ and $f\in C^1(M)$ such that
$\tt\GG(f)(x)=0.$ Applying (\ref{GCD'}) and letting $r\to 0$ we obtain $\GG_2(f)(x)\ge \rr_2\GG^Z(f)(x).$ Assuming further that $\Hess_f(x)=0$ (since one may find such a function $\tt f$ with $\d f(x)=\d\tt f(x)$), it is easy to see from $\GG_2(f)(x)\ge \rr_2\GG^Z(f)(x)$ and $\tt \GG(f)(x)=0$   that $\GG^Z(f)(x)=0.$  Since $\GG+\GG^Z$ is elliptic, we conclude that $\tt \GG(f)(x)=0$ implies $\d f(x)=0$, so that $\tt\GG$ is elliptic as well.
\item[(b)] Even when the family  $\{X_i, \nn_{X_i}X_j\}_{1\le i,j\le n}$ spans the tangent space $TM$,     (\ref{GCD'}) may hold for  some   functions of $r$ in place of   $\rr_1-\ff\kk r$ and $\rr_2$ therein (see Subsection 4.1).
\item[(c)] Combining back to $(\ref{GCD'})$, recall that under condition (\ref{EXX})  it was proved in \cite{BB} (see Proposition 3.1 therein) that
\beq\label{BBG}\GG(P_t f) + \rr_2 t \GG^Z(P_t f)\le \Big(\ff{\rr_2+2\kk}{\rr_2 t} +2\rr_1^-\Big) (P_t f) \big\{P_t (f\log f)-(P_t f)\log P_t f\big\}\end{equation} holds for all $t>0$ and positive $f\in \B_b(M)$.  Comparing with known sharp gradient inequality for the elliptic case, i.e. for $\kk=0$ and $\GG^Z=0$ one has
$$\GG(P_t f) \le \ff{2\rr_1}{\e^{2\rr_1t}-1} (P_t f)\big\{P_t (f\log f)-(P_t f)\log P_t f\big\},$$ where $\ff{2\rr_1}{\e^{2\rr_1t}-1}:=\ff 1 t$ for $\rr_1=0$, the inequality (\ref{BBG}) is less sharp when $\rr_1\ne 0$.   So, it would be nice to find an exact extension of this sharp inequality to the subelliptic setting (see Proposition \ref{C1.2} below).  \end{enumerate}

\

The generalized curvature-dimension condition we proposed here is
$$\GG_2(f)+\sum_{i=1}^lr_i\GG^{(i)}_2(f)\ge \ff{(Lf)^2}d + \sum_{i=0}^l K_i(r_1,\cdots, r_l) \GG^{(i)}(f),\ \ f\in C^3(M), r_1,\cdots, r_l>0,$$ where $d\in (0,\infty]$ is a constant, $\GG^{(0)}:=\GG, \{\GG^{(i)}\}_{1\le i\le l}$ are some $C^2$-diffusion square fields, and $\{K_i\}_{0\le i\le l}$ are some continuous functions on $(0,\infty)^l$. In this paper, we will only consider the   condition with $d=\infty$, i.e.
\beq\label{GC} \GG_2(f)+\sum_{i=1}^lr_i\GG^{(i)}_2(f)\ge   \sum_{i=0}^l K_i(r_1,\cdots, r_l) \GG^{(i)}(f),\ \ f\in C^3(M), r_1,\cdots, r_l>0,\end{equation}
but the condition with finite $d$ will be useful for other purposes as in \cite{BBG,BG1}. In fact,  we will  make use of the following assumption.

\paragraph{(A)}  (\ref{GC})  \emph{ holds for some $C^2$-diffusion square fields $\{\GG^{(i)}\}_{i=0}^l$ and $\{K_i\}_{0\le i\le l}\subset C((0,\infty)^l)$,
where $\GG^{(0)}=\GG$. There exists a smooth compact function $W\ge 1$ on $M$ and a constant
$C>0$ such that $L W\le C W$ and
$\tt\GG (W)\le CW^2$, where $\tt\GG=\sum_{i=0}^l \GG^{(i)}.$  }

\

Recall that  $W$ is called a compact function if $\{W\le r\}$ is compact for any constant $r$. The condition $LW\le CW$ is standard to ensure the non-explosion of the $L$-diffusion process, and the condition $\tt\GG(W)\le C W^2$ is used to prove the boundedness of $\tt\GG(P_t f)$ for $f\in \C$, where
$$\C:= \Big\{f\in C^\infty(M)\cap\B_b(M):\ \ \tt\GG(f)\ \text{is\ bounded}\Big\}.$$
We note that under (\ref{GCD'}) the boundedness of $\tt\GG(P_t f)$ is claimed in \cite{BG1} using a global parabolic comparison theorem, which, however, is not yet available on  general non-compact manifolds.

Similarly to the analysis of elliptic diffusions, a starting point for analyzing the semigroup using a curvature condition is the following $``$gradient" inequalities, which generalize the corresponding ones derived in \cite{BB}. Let

\beg{thm}\label{T1.1} Assume {\bf (A)}. For fixed $t>0$, let $\{b_i\}_{0\le i\le l}\subset C^1([0,t])$ be strictly positive on $(0,t)$ such that
\beq\label{A1} b_i'(s) + 2 \Big\{b_0 K_i\Big(\ff{b_1}{b_0},\cdots,\ff{b_l}{b_0}\Big)\Big\}(s) \ge 0,\ \ s\in (0,t), 1\le i\le l\end{equation} and
$$c_b:= -\inf_{(0,t)} \Big\{b_0'+2 b_0 K_0\Big(\ff{b_1}{b_0},\cdots,\ff{b_l}{b_0}\Big)\Big\}<\infty.$$
Then:
\beg{enumerate} \item[$(1)$] For any $f\in \C$,
$$2\sum_{i=0}^l \big\{b_i(0) \GG^{(i)}(P_tf)-b_i(t)P_t \GG^{(i)}(f) \big\}
\le c_b
\big\{P_tf^2 -(P_t f)^2\big\}.$$
\item[$(2)$] If
\beq\label{EX} \GG^{(i)}\big(\GG(f),f\big)= \GG\big(\GG^{(i)}(f),f\big),\ \ \ 1\le i\le l, f\in C^\infty(M),\end{equation} then for any  positive
$f\in \C$,
$$\sum_{i=0}^l \bigg\{b_i(0) \ff{ \GG^{(i)}(P_tf)}{P_t f}-b_i(t)P_t \ff{\GG^{(i)}(f)}{f}\bigg\}
\le c_b
\big\{P_t(f\log f)-(P_t f)\log P_t f\big\}.$$ \end{enumerate}\end{thm}

Theorem \ref{T1.1}  will be proved in Section 2. In Section 3 we  apply this theorem to the study of $``$gradient" estimate, Poincar\'e inequality, Harnack inequality and applications.
In Section 4 we present some specific examples to illustrate the generalized curvature condition (\ref{GC}), for which (\ref{GCD'}) does not hold. We will not consider the log-Sobolev and HWI inequalities for    enlarged Dirichlet forms given by $\sum_{i=0}^l\GG^{(i)}$, since they are no longer intrinsic for the underlying subelliptic diffusion process.  It is still open to derive the intrinsic semigroup log-Sobolev inequality and HWI inequality for subelliptic diffusion processes using generalized curvature conditions.

Finally, noting that (\ref{GC}) is still not available for some highly degenerate subelliptic diffusion operators, e.g. $L=\ff{\pp^2}{\pp x^2}+x\ff{\pp}{\pp y}$
on $\R^2$ for which the derivative formula and Harnack ineuqlity has been established in \cite{Z,GW}, we propose in Section 5 an extension of Theorem \ref{T1.1} with (\ref{GC}) holding for not necessarily positively definite $\{\GG^{(i)}\}_{1\le i\le l}$ and not necessarily all $r_i>0.$ As application, explicit
gradient-entropy inequality and Harnack inequality are presented for this simple example.

\section{Proof of Theorem 1.1}

To prove this theorem using a modified Bakry-Emery semigroup argument as in \cite{BB}, we need to first confirm that $P_t\C\subset \C$, which follows immediately from the following lemma.

\beg{lem}\label{L2.1} Assume {\bf (A)} and let $K=\min_{0\le i\le l} K_i(1,\cdots,1).$ Then
\beq\label{GE} \tt\GG(P_t f)\le \e^{-2Kt} P_t\tt\GG(f),\ \ t\ge 0, f\in \C.\end{equation}
\end{lem}

\beg{proof} (i) We first prove for any $f\in C_0^2(M)$ and $t>0$, $\tt\GG(P_\cdot f)$ is bounded on $[0,t]\times M.$ To this end, we approximate the generator $L$ by using operators with compact support, so that the approximating diffusion processes stay in  compact sets. Take $h\in C_0^\infty([0,\infty))$ such that $h'\le 0, h|_{[0,1]}=1$ and supp$h=[0,2].$ For any $m\ge 1,$ let $\varphi_m = h(W/m)$ and
$L_m= \varphi_m^2 L.$ Then $L_m$ has compact support $B_m:=\{W\le 2m\}.$  Let $x\in \{W\le m\}$ and $X_s^m$ be the $L_m$-diffusion process starting at $x$. Let
$$\tau_m= \inf\{s\ge 0: W(X_s^m)\ge 2m\}.$$ Since $LW\le CW, \GG(W)\le \tt\GG(W)\le CW^2, h'\le 0, 0\le h\le 1$ and $h'(W/m)=0$ for $ W\ge 2m$, we have
\beg{equation*}\beg{split}&L_m\ff 1 {\varphi_m^2} =-\ff{2L\varphi_m}{\varphi_m} + \ff{6\GG(\varphi_m)}{\varphi_m^2}\\
& = -\ff{2h'(W/m)LW}{m\varphi_m}-\ff{2h''(W/m)\GG(W)}{m^2\varphi_m} +\ff{6h'(W/m)^2\GG(W)}{m^2\varphi_m^2}
\le \ff{C_1}{\varphi_m^2}\end{split}\end{equation*} for some constant $C_1>0$ independent of $m$. By a standard argument, this implies that $\tau_m=\infty$ and
\beq\label{NE}\E\Big(\ff 1 {\varphi_m^2} \Big)(X_s) \le \ff{\e^{C_1 t}}{\varphi_m^2(x)}=\e^{C_1 s}, \ \ \ s\ge 0.\end{equation} Now, let $P_s^m$ be the diffusion semigroup generated by $L_m$. By the It\^o formula and $\tt\GG_2\ge K\tt\GG$ implied by {\bf (A)} we obtain
\beg{equation}\label{EEG}\beg{split}\d \tt\GG(P_{t-s}^mf)(X_s^m) &=\d M_s^m + \big\{\varphi_m^2 L\tt\GG(P_{t-s}^mf)- 2\tt\GG(P_{t-s}^mf, \varphi_m^2 LP_{t-s}^mf)\big\}(X_s^m)\d s\\
&\ge \d M_s^m+ \big\{2 \varphi_m^2\tt\GG_2(P_{t-s}^mf) -4\tt\GG(\log \varphi_m, P_{t-s}^mf) P_{t-s}^m L_m f\big\}(X_s^m)\d s \\
&\ge \d M_s^m - \Big\{2|K|\tt\GG(P_{t-s}^mf)+4\|Lf\|_\infty \ss{\tt\GG(\log \varphi_m)\tt\GG(P_{t-s}^mf)}\Big\}(X_s^m)\d s\\
&\ge \d M_s^m -C_2 \tt\GG(P_{t-s}^mf)(X_s^m)\d s - \tt\GG(\log\varphi_m)(X_s^m)\d s,\ \ s\in [0,t] \end{split}\end{equation} for some martingale $M_s^m$ and some constant $C_2>0$ independent of $m$. Since $h'(W/m)=0$ for $W\ge 2m$ and $\tt\GG(W)\le CW^2$,   $$\tt\GG(\log\varphi_m)=\ff{h'(W/m)^2\tt\GG(W)}{m^2\varphi_m^2}\le \ff{C_3}{\varphi_m^2}$$ holds for some constant $C_3>0$ independent of $m$. Combining this with (\ref{NE}) and (\ref{EEG}) we conclude that
\beq\label{UG}\tt\GG(P_t^mf)\le \e^{C_2 t}P_t^m\tt\GG(f) + C_3\int_0^t \E\Big(\ff 1 {\varphi_m^2}\Big)(X_s^m)\d s \le \e^{C_2t}\|\tt\GG(f)\|_\infty+C_3\e^{C_1t} \end{equation} holds on $\{W\le m\}.$
Letting $\tt\rr$ be the intrinsic distance induced by $\tt\GG$, i.e.
$$\tt\rr(z,y):= \sup\{|g(x)-g(y)|:\ \tt\GG(g)\le 1\},\ \ z,y\in M,$$ we deduce from (\ref{UG})   that for any $x,y\in M$,
\beq\label{LL} |P_t^mf(z)-P_t^mf(y)|^2\le \tt\rr(z,y)^2 \Big(\e^{C_2t}\|\tt\GG(f)\|_\infty+C_3\e^{C_1t}\Big),\ \ t> 0\end{equation} holds for large enough $m$. Noting that the $L$-diffusion process is non-explosive and $X_s^m$ is indeed generated by $L$ before time $\si_m:=\inf\{s\ge 0: W(X_s^m)\ge m\}$ which increases  to $\infty$ as
$m\to\infty$,  we conclude that $\lim_{m\to\infty}P_t^mf= P_tf$ holds point-wisely. Therefore, letting $m\to\infty$ in (\ref{LL}) we obtain
$$|P_tf(z)-P_tf(y)|^2\le \tt\rr(z,y)^2 \Big(\e^{C_2t}\|\tt\GG(f)\|_\infty+C_3\e^{C_1t}\Big),\ \ t\ge 0,y,z\in M.$$ This implies that $\tt\GG(P_\cdot f)$ is bounded on $[0,t]\times M$ for any $t>0.$

(ii) By an approximation argument, it suffices to prove (\ref{GE}) for $f\in C_0^2(M)$. By the It\^o formula and (\ref{GC}), there exists a local martingale $M_s$ such that
 $$ \d \tt\GG(P_{t-s} f)(X_s) =\d M_s + 2\tt\GG_2(P_{t-s} f)(X_s)\d s
 \ge \d M_s + 2 K \tt\GG(P_{t-s}f)(X_s)\d s,\ \ s\in [0,t].$$ Thus,
 $$[0,t]\ni s\mapsto \e^{-2 Ks}\tt\GG(P_{t-s}f)(X_s)$$ is a local submartingale. Since due to (a) this process is bounded, so that it is indeed a submartingale. Therefore, (\ref{GE}) holds.
\end{proof}

%By an approximation argument, we only need to prove (\ref{GE}) for $f\in C^\infty_0(M).$ Indeed, let
%$$\tt\rr(x,y)= \sup\big\{|f(x)-f(y)|:\ f\in C^1(M), \tt\GG(f)\le 1\big\},\ \ x,y\in M.$$Then (\ref{GE}) is equivalent to
%\beq\label{GE'} \sup_{U_x\setminus\{x\}} \ff{|P_t f(x)-P_t f|^2}{\rr(x,\cdot)^2} \le \e^{-2Kt} \sup_{U_x} P_t \tt\GG(f),\ \ x\in M, U_x \text{\ is\ any\ neighborhood\ of\ }x.\end{equation} Now,  for any $f\in C^1(M)\cap\B_b(M)$ such that $ \|\tt\GG(f)\|_\infty<\infty$, there exists a sequence $\{f_m\}_{m\ge 1}\subset C_0^\infty(M)$ such that $f_m\to f$ uniformly, $\tt\GG(f_m)\to\tt\GG(f)$ locally uniformly, and $\sup_{m\ge 1}\|\tt\GG(f_m)\|_\infty<\infty.$
% Then   (\ref{GE'}) holds provided it holds for all $f_m$ in place of $f$.
%Now, for $f\in C_0^\infty(M)$, one has $\ff{\d}{\d t} P_t f= LP_t f, \ff{\d}{\d t} P_t\tt\GG( f)= LP_t \tt\GG(f),$ so that due to (\ref{C}),
% $$\varphi_t:= \e^{2Kt} \tt\GG(P_t f),\ \ \psi_t:= P_t\tt\GG(f)$$ satisfy
% \beg{equation*}\beg{split} &\ff{\d}{\d t} \psi_t = LP_t \tt\GG(f)= L\psi_t,\\
% &\ff{\d}{\d t} \varphi_t = L\varphi_t-2\e^{2Kt}\tt\GG_2(P_tf)+2K\e^{2Kt}\tt\GG(P_t f)\le L\varphi_t.\end{split}\end{equation*} Therefore, $\Phi:=\varphi-\psi$ satisfies $\pp_t \Phi \le L\Phi, \Phi_0=0.$ Since the $L$-diffusion process is non-explosive so that the heat equation associated to $L$ has a unique solution, it follows from the comparison theorem for parabolic partial differential equations that
 %$\Phi_t\le 0$, i.e. $\varphi_1(t)\le \varphi_2(t)$  (cf. page 52 in \cite{F}).

\beg{proof}[Proof of Theorem \ref{T1.1}] (1)   It suffices to prove for $f\in C^\infty(M)$ which is constant outside a compact set.
In this case we have $\ff{\d}{\d s}P_s f= LP_s f= P_s Lf.$ Since $X_s$ is non-explosive, by the It\^o formula for any $0\le i\le l$ there exists a local martingale $M_s^{(i)}$ such that
\beg{equation*}\beg{split} \d\GG^{(i)} (P_{t-s}f)(X_s) &= \d M_s^{(i)} +\big\{L\GG^{(i)}(P_{t-s}f)-2 \GG^{(i)}(P_{t-s}f, LP_{t-s}f)\big\}(X_s)\d s\\
&= \d M_s^{(i)} + 2 \GG_2^{(i)}(P_{t-s}f)(X_s)\d s,\ \ \ s\in [0,t].\end{split}\end{equation*}
Therefore, due to (\ref{GC}) and (\ref{A1}), there exists a local martingale $M_s$ such that
\beg{equation*}\beg{split} &\d\bigg\{\sum_{i=0}^l b_i(s) \GG^{(i)}(P_{t-s}f)(X_s)\bigg\}\\
&\ge \d M_s + \bigg\{\sum_{i=0}^l \Big(2b_i(s)\GG_2^{(i)}(P_{t-s}f)+b_i'(s)\GG^{(i)}(P_{t-s}f)\Big)\bigg\}(X_s)\d s\\
&\ge \d M_s +\bigg\{ \sum_{i=0}^l \Big(b_i'(s) +2b_0(s)K_i\Big(\ff{b_1}{b_0},\cdots,\ff{b_l}{b_0}\Big)(s)\Big)\GG^{(i)}(P_{t-s}f)\bigg\}(X_s)\d s\\
&\ge \d M_s +\bigg\{b_0'(s) +2b_0(s)K_0\Big(\ff{b_1}{b_0},\cdots,\ff{b_l}{b_0}\Big)(s)\bigg\} \GG((P_{t-s}f) (X_s)\d s.\end{split}\end{equation*}
So, if $c_b<\infty$
then
$$\sum_{i=0}^l b_i(s)   \GG^{(i)}(P_{t-s}f) (X_s)+c_b\int_0^s  \GG(P_{t-r}f) (X_r)\d r$$ is a local submartingale for $s\in [0,t].$
Since, due to Lemma \ref{L2.1}, $\{\GG^{(i)}(P_{t-s}f)\}_{0\le i\le l}$ are bounded, it is indeed a   submartingale. In particular,
$$\sum_{i=0}^l \big\{b_i(0) \GG^{(i)}(P_tf)-b_i(t)P_t \GG^{(i)}(f) \big\}\le c_b \int_0^t P_s \GG(P_{t-s}f)  \d s.$$
Then the proof is finished by noting that
$$P_s  \GG(P_{t-s}f) =  \ff 1 2  \ff{\d}{\d s} P_s (P_{t-s}f)^2.$$

(2) Let $f$ be strictly positive and be constant outside a compact set. Let
$$\phi^{(i)}(s,x)= \big\{(P_{t-s }f) \GG^{(i)} (\log P_{t-s}f)\big\}(x),\ \ \ 0\le i\le l, s\in [0,t], x\in M.$$ It is easy to see that (\ref{EX}) implies (cf. \cite{BG1})
$$L\phi^{(i)} +\ff{\pp}{\pp s} \phi^{(i)} =2 (P_{t-s}f) \GG_2^{(i)}(\log P_{t-s} f),\ \ 0\le i\le l.$$
So,  for each $0\le i\le l,$ there exists a local martingale $M_s^{(i)}$ such that
$$\d\phi^{(i)}(s,X_s)
 = \d M_s^{(i)} +  2 \big\{(P_{t-s}f) \GG_2^{(i)}(\log P_{t-s} f)\big\}(X_s)\d s,\ \ s\in [0,t].$$  The remainder of the proof is then completely similar to (1); that is,
 $$\sum_{i=0}^l b_i(s) \big\{(P_{t-s}f)  \GG^{(i)}(\log P_{t-s}f)\big\}(X_s)+c_b\int_0^s \big\{(P_{t-r}f)\GG(\log P_{t-r}f))\big\}(X_r)\d r$$ is a   submartingale for $s\in [0,t]$, so that the desired inequality follows by noting that
$$P_s \big\{(P_{t-s}f)\GG(\log P_{t-s}f)\big\}=   \ff{\d}{\d s} P_s\big\{(P_{t-s}f) \log P_{t-s}f\big\}.$$\end{proof}

\section{Applications of Theorem \ref{T1.1}}

For any non-negative symmetric measurable functions $\tt\rr$ on $M\times M$, let $W_2^{\tt\rr}$ be the  $L^2$-transportation-cost with cost function
$\tt\rr;$ i.e. for any two probability measures $\mu_1,\mu_2$ on $M$,
$$W_2^{\tt\rr}(\mu_1,\mu_2):= \inf_{\pi\in\C(\mu_1,\mu_2)} \pi(\tt\rr^2)^{1/2},$$ where $\pi(\tt\rr)$ stands for the integral of $\tt\rr$ w.r.t. $\pi$, and
$\C(\mu_1,\mu_2)$ is the set of all couplings of $\mu_1$ and $\mu_2$.

\subsection{$L^2$-derivative estimate and applications}

\beg{prp}\label{P3.1} Assume {\bf (A)}. Let $t>0$ and $\{b_i\}_{0\le i\le l}\subset C^1([0,t])$ be strictly positive in $(0,t)$ such that $(\ref{A1})$ holds. If
$b_i(t)=0, 0\le i\le l$ and $c_b<\infty,$  then:
\beg{enumerate} \item[$(1)$] For any $f\in \B_b(M),$
\beq\label{3.1} \sum_{i=0}^l b_i(0)\GG^{(i)} (P_t f) \le c_b \big\{P_t f^2-(P_tf)^2\big\}.\end{equation}
\item[$(2)$] For any non-negative $f\in \B_b(M)$, the Harnack type inequality
\beq\label{3.2} P_t f(x)\le P_t f(y) +c_b \rr_b(x,y) \ss{P_t f^2(x)},\ \ x,y\in M\end{equation} holds for $\rr_b$ being the intrinsic distance induced by
$\GG_b:= \sum_{i=0}^l b_i(0) \GG^{(i)}.$
\item[$(3)$] If $P_t$ has an   invariant probability  measure $\mu$, then for any $f\ge 0$ with $\mu(f)=1$, the variance-cost inequality
\beq\label{VC}{\rm Var}_{\mu}(P_t^* f)\le c_b W_2^{\rr_b}(f\mu,\mu)\ss{\mu((P_t^*f)^3)} \end{equation} holds, where $P_t^*$ is the adjoint operator of $P_t$ in $L^2(\mu)$, and
$${\rm Var}_{\mu}(P_t^* f):= \mu((P_t^*f)^2)-\mu(P_t^*f)^2= \mu((P_t^*f)^2)-1.$$\end{enumerate}\end{prp}

\beg{proof} The first assertion is a direct consequence of Theorem \ref{T1.1}(1), while according to \cite[Proposition 1.3]{W11} (\ref{3.1}) is equivalent to (\ref{3.2}). Finally, (\ref{VC})  follows from (\ref{3.2}) according to  the following lemma \ref{L3.2}.\end{proof}

\beg{lem}\label{L3.2} Let $P$ be a Markov operator on $\B_b(E)$ for a measurable space $(E,\B)$. Let $\mu$ be an invariant probability measure of $P$. If \beq\label{HH} Pf(x)\le Pf(y) +C \rr(x,y) \ss{Pf^2(x)},\ \ f\in \B_b^+(E)\end{equation} holds for some constant $C>0$ and non-negative symmetric function $\rr$ on $E\times E$, then
$$ {\rm Var}_\mu(P^*f)\le  C    W_2^{\rr} (f\mu,\mu)  \ss{ \mu((P^*f)^3)},\ \ f\ge 0, \mu(f)=1.$$\end{lem}

\beg{proof} Let $f\ge 0$ with $\mu(f)=1$. For any $\pi\in \C(f\mu,\mu)$, (\ref{HH}) implies
\beg{equation*}\beg{split} \mu((P^*f)^2) &= \mu(fPP^*f)=\int_{E\times E} P(P^*f)(x)\pi(\d x,\d y)\\
& \le \int_{E\times E} P(P^*f)(y) \pi(\d x,\d y)
+C\int_{E\times E} \rr(x,y)\ss{P(P^*f)^2(x)}\,\pi(\d x,\d y)\\
&\le \mu(PP^*f) +C\ss{\pi(\rr^2)\mu(fP(P^*f)^2)}= 1 + C\ss{\pi(\rr^2) \mu((P^*f)^3)}.\end{split}\end{equation*} This   completes the proof.
\end{proof}

\subsection{Entropy-derivative estimate and applications}
\beg{prp}\label{P3.12} Assume {\bf (A)} and $(\ref{EX})$. Let $t>0$ and $\{b_i\}_{0\le i\le l}\subset C^1([0,t])$ be strictly positive in $(0,t)$ such that $(\ref{A1})$ holds. If
$b_i(t)=0, 0\le i\le l$ and
$c_b<\infty,$ then:
\beg{enumerate} \item[$(1)$] For any $f\in \B_b(M),$
\beq\label{3.1'} \sum_{i=0}^l b_i(0)\GG^{(i)} (P_t f) \le c_b (P_tf)\big\{P_t(f\log f)- (P_tf)\log P_t f\big\}.\end{equation}
\item[$(2)$] For any non-negative $f\in \B_b(M)$  and $\aa>1$, the Harnack type inequality
\beq\label{3.2'} (P_t f)^\aa (x)\le P_t f^\aa (y) \exp\bigg[\ff{\aa c_b\rr_b(x,y)^2}{4(\aa-1)}\bigg],\ \ x,y\in M\end{equation} holds for $\rr_b$ being the intrinsic distance induced by
$\GG_b:= \sum_{i=0}^l b_i(0) \GG^{(i)}.$ Consequently, the log-Harnack inequality
\beq\label{LH} P_t\log f(x)\le \log P_t f(y) + \ff{c_b\rr_b(x,y)^2}{4}\end{equation} holds for uniformly positive measurable function $f.$
\item[$(3)$] If $P_t$ has an   invariant probability  measure $\mu$, then for any $f\ge 0$ with $\mu(f)=1$, the entropy-cost inequality
\beq\label{EC} \mu \big((P_t^* f)\log P_t^*f\big) \le \ff  {c_b^2}4     W_2^{\rr_b} (f\mu,\mu)^2.\end{equation} \end{enumerate}\end{prp}

\beg{proof} The first assertion is a direct consequence of Theorem \ref{T1.1}(2), (\ref{3.2'}) follows from (1) and the following Lemma \ref{L3.4}, (\ref{LH}) follows from (\ref{3.2'}) according to \cite[Proposition 2.2]{W10b}, and finally, (\ref{EC}) follows from (\ref{LH}) according to \cite[Proposition 4.6]{W11b}.\end{proof}

\beg{lem}\label{L3.4} Let $P$ be a Markov operator on $\B_b(M)$ and let $\gg$ be a positive measurable function on $(0,\infty)$. Let $\bar\GG$ be a smooth diffusion square field on $E$ with intrinsic distance $\bar\rr$. If \beq\label{G} \ss{\bar\GG(Pf)}\le \dd \big\{P(f\log f)-(Pf)\log Pf\big\}+\gg(\dd) Pf,\ \ \dd>0 \end{equation} holds for  all positive $f\in \B_b(M)$, then for any $\aa>1$,
$$(Pf)^\aa(x)\le \big(Pf^\aa(y)\big)\exp\bigg[\int_0^1  \ff{\aa\bar\rr(x,y)}{1+(\aa-1)s} \gg\Big(\ff{\aa-1}{(1+(\aa-1)s)\bar\rr(x,y)}\Big)\d s\bigg],\ \ x,y\in M$$ holds for all positive $f\in \B_b(M).$\end{lem}

\beg{proof} The proof is completely similar to that of \cite[Theorem 1.2]{ATW09}. Let $\bar\rr(x,y)<\infty$, and let $x_\cdot: [0,1]\to M$ with constant speed w.r.t. $\bar\rr$ such that $x_0=x, x_1=y$; that is, the curve is  the minimal geodesic from $x$ to $y$ induced by the metric associated to $\tt\GG$.
Then  $|\ff{\d f(x_s)}{\d s}|^2\le \bar\rr(x,y)^2\bar\GG(f)(x_s)$ holds for all $f\in C^1(M)$.  Let $\bb(s)=1+(\aa-1)s.$ It follows from (\ref{G}) that
\beg{equation*}\beg{split} &\ff{\d}{\d s} \log (Pf^{\bb(s)})^{\ff\aa{\bb(s)}}(x_s)\\
 &\ge \ff{\aa(\aa-1)\{P(f^{\bb(s)}\log f^{\bb(s)}) -(Pf^{\bb(s)})\log Pf^{\bb(s)}\}}{\bb(s)^2Pf^{\bb(s)}}(x_s) -\ff{\aa \bar\rr(x,y)\ss{\bar \GG(f)}}{\bb(s)Pf^{\bb(s)}}(x_s)\\
&\ge - \ff{\aa\bar\rr(x,y)}{\bb(s)} \gg\Big(\ff{\aa-1}{\bb(s)\bar\rr(x,y)}\Big).\end{split}\end{equation*} Then the proof is finished by integrating both sides on $[0,1]$ w.r.t. $\d s$.\end{proof}
For applications of the Harnack and log-Harnack inequalities, we   refer to Section 4 of \cite{W11b} (see also \cite{W10b, WY}).  In particular, if $P_t$ is symmetric w.r.t. some probability measure $\mu$ such that $\mu(\e^{\ll\rr_b(o,\cdot)^2})<\infty$ holds for some $\ll>\ff {c_b}4,$ then (\ref{3.2'}) implies that the log-Sobolev inequality
$$\mu(f^2\log f^2)\le c \mu(\GG(f)),\ \ f\in C_0^1(M),\ \mu(f^2)=1$$ holds for some constant $c>0.$

\subsection{Exponential decay and  Poincar\'e inequality}

\beg{prp} \label{P3.5}Assume {\bf (A)}. For   $r_i>0, 1\le i\le l$, let $$\ll(r_1,\cdots, r_l)=\min_{0\le i\le l} \ff{K_i(r_1,\cdots, r_l)}{r_i},$$ where $r_0:=1.$
Then
$$\sum_{i=0}^l r_i\GG^{(i)} (P_t f) \le \e^{-2\ll(r_1,\cdots, r_l) t}\sum_{i=0}^l r_i  P_t \GG^{(i)}(f),\ \ t\ge 0, f\in C_b^1(M).$$
Consequently, if $P_t$ is symmetric w.r.t a probability measure $\mu$ and
$$ \ll:=\sup_{r_1,\cdots, r_l>0} \ll(r_1,\cdots r_l)>0,$$ then the Poinca\'re inequality
\beq\label{PP}\mu(f^2)\le \ff 1 \ll \mu(\GG(f))+\mu(f)^2,\  \ f\in C_0^1(M)\end{equation} holds. \end{prp}
\beg{proof} By a standard spectral theory (cf. the proof of \cite[Corollary 2.4]{BB}), the Poincar\'e inequality follows immediately from the desired derivative inequality. To prove the derivative inequality, we take
$$b_0(s)= \e^{-2\ll(r_1,\cdots, r_l) s},\ \ b_i(s)= r_i b_0(s), \ \ \ 1\le i\le l, s\ge 0.$$ Then
$$b_i' + 2 b_0K_i\Big(\ff{b_1}{b_0},\cdots, \ff{b_l}{b_0}\Big) = -2r_i \ll(r_1,\cdots, r_l) b_0 + 2 b_0K_i(r_1,\cdots, r_l)\ge 0,\ \ 0\le i\le l.$$ Therefor, the desired gradient inequality follows from Theorem \ref{T1.1}(1). \end{proof}

When (\ref{GCD'}) holds for $\rr_1,\rr_2>0$ and $\kk\ge 0$, we have
$$\ll = \sup_{r>0} \Big\{\ff {\rr_2} r\land \Big(\rr_1-\ff\kk r\Big)\Big\}= \ff{\rr_1\rr_2}{\rr_2+\kk}>0.$$ Thus, (\ref{PP})
recovers the Poincar\'e inequality presented in \cite[Corollary 2.4]{BB}.

\subsection{Refined derivative inequalities under (\ref{GCD'})}
Coming back to condition (\ref{GCD'}),   Theorem \ref{T1.1} implies the following exact extensions of sharp gradient estimates in the elliptic setting (see the above comment (c)).

\beg{prp}\label{C1.2} Assume $(\ref{GCD'})$ for some constants $\rr_2>0,\kk\ge 0$ and $\rr_1\in\R$. Assume there exist a smooth compact function $W\ge 1$ and a constant $C>0$ such that $LW\le CW$ and $\tt\GG(W)\le CW^2$.

\beg{enumerate} \item[$(1)$] For any $t>0$ and $f\in \B_b(M)$,
\beg{equation*}\beg{split} &\GG(P_t f)+ \ff{\rr_2(\e^{2\rr_1 t}-1-2\rr_1t)}{\rr_1(\e^{2\rr_1 t}-1)} \GG^Z(P_t f)\\
&\le \Big(1+ \ff{\kk (\e^{2\rr_1^+t}-1)^2}{\rr_2(\e^{2\rr_1^+t}-1-2\rr_1^+t)}\Big)\ff{\rr_1}{\e^{2\rr_1 t}-1}\big\{P_t f^2-(P_t f)^2\big\},\end{split}\end{equation*} where when $\rr_1\le 0$,
$$\ff{(\e^{2\rr_1^+t}-1)^2}{\e^{2\rr_1^+t}-1-2\rr_1^+t}:=\lim_{r\downarrow 0} \ff{(\e^{r}-1)^2}{\e^{r}-1-r}=2.$$ Consequently, if $\rr_1>0$ and
$P_t$ is symmetric w.r.t. a probability measure $\mu$, then the Poincar\'e inequality
\beq\label{P}\mu(f^2)\le \ff 1 {\rr_1} \mu(\GG(f)) +\mu(f)^2,\ \ f\in C_0^1(M)\end{equation} holds.
\item[$(2)$] If $(\ref{EX})$ holds, then for any $t>0$ and positive $f\in \B_b(M)$,
\beg{equation*}\beg{split} &\GG(P_t f)+ \ff{\rr_2(\e^{2\rr_1 t}-1-2\rr_1t)}{\rr_1(\e^{2\rr_1 t}-1)} \GG^Z(P_t f)\\
&\le \Big(1+ \ff{\kk (\e^{2\rr_1^+t}-1)^2}{\rr_2(\e^{2\rr_1^+t}-1-2\rr_1^+t)}\Big)\ff{2\rr_1(P_tf) \{P_t (f\log f) -(P_t f)\log P_t f\}}{\e^{2\rr_1 t}-1}.\end{split}\end{equation*} \end{enumerate}\end{prp}

\beg{proof}  Let
\beg{equation*}\beg{split} &b_0(s)= \ff{\e^{2\rr_1(t-s)}-1}{2\rr_1},\\
&b_1(s) =2\rr_2\int_s^t b_0(r)\d r = \ff{\rr_2(\e^{2\rr_1(t-s)}-1-2\rr_1(t-s)}{2\rr_1^2},\ \ s\in [0,t].\end{split}\end{equation*}
Then it is easy to see that $b_1+ 2b_0\rr_2=0$ and
\beg{equation*}\beg{split} \Big\{b_0'+2b_0\Big(\rr_1-\ff{Kb_0}{b_1}\Big)\Big\}(s)&= -1-\ff{\kk(\e^{2\rr_1(t-s)}-1)^2}{\rr_2(\e^{2\rr_1(t-s)}-1-2\rr_1(t-s))}\\
&\ge -1-\ff{\kk (\e^{2\rr_1^+t}-1)^2}{\rr_2(\e^{2\rr_1^+t}-1-2\rr_1^+t)}.\end{split}\end{equation*} Since  (\ref{GCD'}) implies (\ref{GC}) for
$l=1, \GG^{(1)}=\GG^Z, K_0(r)=\rr_1-\ff{\kk} r $ and $K_1(r)=\rr_2,$ the desired derivative inequalities follows from (\ref{3.1}) and (\ref{3.1'}).\end{proof}

\section{Examples} Additionally to  those examples given in \cite{BG1} such that (\ref{GCD}) holds, we   present three simple  examples for which (\ref{GCD'})\,(and hence (\ref{GCD})) is not available  but   our more general condition (\ref{GC}) holds true.  For the first example the Poincar\'e and log-Sobolev inequalities are confirmed in the symmetric setting, while for the other two examples (\ref{EX}) does not hold so that we are only able to derive results in Proposition \ref{P3.1}.

\subsection{Example A}  Let $M=\R\times \bar M$, where $\bar M$ is a complete connected Riemannian manifold. Let $\bar L$ be an elliptic differential operator on $\bar M$ satisfying the curvature-dimension condition
\beq\label{CD1} \bar\GG_2(f)\ge K\bar\GG(f)+\ff{(\bar Lf)^2}{m},\ \ f\in C^\infty(\bar M),\end{equation} for some constant $K\ge 0$ and $m\in (1,\infty)$, where $\bar\GG$ is the square field of $\bar L$ and $\bar\GG_2$ is the associated curvature operator, i.e.
$\bar\GG_2(f)= \ff 1 2 \bar L\bar\GG(f)- \bar\GG(f, Lf).$  Consider
$$ Lf (x,y)=   f_{xx}(x,y) - r_0 x  f_x(x,y) +x^2 \bar Lf(x,\cdot)(y),\ \ f\in C^\infty(M), (x,y)\in M$$ for some constant $r_0\in\R,$ where and in the sequel, $f_{x_1\cdots x_k}:=\ff{\pp^k}{\pp x_1\cdots\pp x_k}f,\ k\ge 1.$ Then
$$\GG^{(0)}(f,g)(x,y):=\GG(f,g)(x,y)= (f_xg_x)(x,y)+ x^2\bar\GG(f(x,\cdot),g(x,\cdot))(y).$$
Let
$$\GG^{(1)}(f,g)(x,y)= \bar\GG(f(x,\cdot),g(x,\cdot))(y),\ \ f\in C^\infty(M), (x,y)\in M.$$
According to (\ref{CD1}), there exists a positive smooth compact function $\bar W$ on $\bar M$ such that $\bar L W, \bar\GG(W)\le 1.$ In fact, let $\bar\rr$ be the intrinsic distance to a fixed point induced by $\bar\GG$, by (\ref{CD1}) for $K\ge 0$ and    the comparison theorem, one has (see \cite{Qian})
$$\bar L\bar\rr\le \ff {m-1}{\bar\rr}$$ outside the fixed point and the cut-locus of this point. By Greene-Wu's  approximation theorem (see \cite{GrW}), we may assume that $\bar\rr^2$ is smooth so that
$\bar L \ss{1+\bar\rr^2}\le c_1$ holds for some constant $c_1>0$. Noting that $\bar\GG(\bar\rr)=1$, we may take $\bar W=\vv \ss{1+\bar\rr^2}$ for small enough constant $\vv>0.$

Now, let  $W(x,y)= 1+x^2+\bar W(y)$, which is a smooth compact function on $M$. It is easy to see that
\beq\label{W1}   LW(x,y)\le   2(1+r_0^-)W(x,y),\ \  \tt\GG(W)(x,y)=4 x^2 +(1+x^2)\bar\GG(\bar W)(y)\le 5W(x,y),\end{equation}
where $\tt\GG= \GG+\GG^{(1)}.$
\beg{prp}\label{P4.1} The generalized curvature condition $(\ref{GC})$ holds for $l=1$ and
$$K_1(r)=1,\ \ \ K_0(r)= \Big(r_0-\ff m r\Big)\land \Big(Kr-r_0-\ff 4 r\Big),\ \ r>0,$$ and $(\ref{EX})$ holds. Consequently:
\beg{enumerate}
\item[$(1)$] Propositions \ref{P3.1} and \ref{P3.12} hold for $b_0(0)=t, b_1(0)=t^2$ and
$$c_b= 1+2\sup_{r\in (0,t)}\big\{(m-r_0r\big)\lor (r_0r+4-Kr^2)\big\}.$$
 \item[$(2)$] If $K,r_0>0$ and $\bar L$ is symmetric w.r.t. a probability measure $\bar\mu$ on $\bar M$, then $P_t$ is symmetric w.r.t.
$$\mu(\d x,\d y):= \Big(\ff{\ss {r_0}\exp[-\ff{r_0}2 x^2]}{\ss{2\pi}}\,\d x\Big)\bar\mu(\d y),$$ and the Poincar\'e inequality $(\ref{PP})$ holds for
$$\ll= \min\Big\{\ff{2K}{r_0+\ss{r_0^2+20K}},\ \ff{r_0}{m+1}\Big\}>0.$$ Moreover,  the log-Sobolve inequality
$$\mu(f^2\log f^2)\le c\mu(\GG(f)),\ \ f\in C_0^1(M), \mu(f^2)=1$$ holds for some constant $c>0.$ \end{enumerate}\end{prp}

\beg{proof} (i) The proof of (\ref{EX}) is trivial. Below we intend to  prove (\ref{GC}) for the desired $K_0$ and $K_1$; that is,
\beq\label{GGG} \GG_2(f)+r \GG^{(1)}_2(f) \ge \GG^{(1)}(f) + \Big\{\Big(r_0-\ff m r\Big)\land \Big(Kr-r_0-\ff 4 r\Big)\Big\}\GG(f),\ \ f\in C^\infty(M).\end{equation}  It is easy to see that at point $(x,y),$
\beg{equation*}\beg{split} &\GG_2(f) = f_{xx}^2  +(1-r_0x^2)\GG^{(1)}(f)
 +4x \GG^{(1)}(f,f_x) \\
 &\qquad\qquad   + 2 x^2 \GG^{(1)}(f_x)  + x^4 \bar \GG_2(f(x,\cdot))(y) -2 x f_x\bar L f(x,\cdot)(y)  + r_0 f_x^2,\\
&\GG_2^{(1)}(f) = x^2 \bar\GG_2(f(x,\cdot) (y) +\GG^{(1)}(f_x).\end{split}\end{equation*} Combining these with (\ref{CD1})we obtain
\beg{equation*}\beg{split} & \GG_2(f)+r \GG^{(1)}_2(f)\\
&\ge \GG^{(1)}(f)-r_0x^2\GG^{(1)}(f) +\big\{(2x^2+r)\GG^{(1)}(f_x)+4 x \GG^{(1)}(f,f_x)\big\} + (x^4+rx^2)K\GG^{(1)}(f)\\
&\quad +\bigg\{\ff{(x^4+rx^2)(\bar Lf(x,\cdot)(y))^2}m -2 x f_x \bar Lf(x,\cdot)(y)\bigg\}+r_0 f_x^2\\
&\ge \GG^{(1)}(f) + \Big(K(x^2+r)-r_0-\ff 4 {2 x^2 +r}\Big)x^2\GG^{(1)}(f) -\ff{mx^2}{x^4+rx^2}f_x^2+r_0f_x^2\\
&\ge \GG^{(1)}(f) + \Big(Kr-r_0-\ff 4 r\Big)x^2 \GG^{(1)}(f) +\Big(r_0- \ff m r\Big)f_x^2\\
&\ge \GG^{(1)}(f) + \Big\{\Big(r_0-\ff m r\Big)\land \Big(Kr-r_0-\ff 4 r\Big)\Big\}\GG(f).\end{split}\end{equation*} Therefore, (\ref{GGG}) holds.

(ii) Whence (\ref {EX}) and (\ref{GC}) are confirmed for the desired $K_0$ and $K_1$, due to (\ref{W1}) the assumption {\bf (A)} holds. Then (1) follows immediately  by taking
$$b_0(s)= t-s,\ \ b_1(s)=(t-s)^2,\ \ s\in [0,t].$$ It remains to prove the  the Poincar\'e inequality and the log-Sobolev inequality for $K,r_0>0$ in the symmetric setting. By Proposition \ref{P3.5}, the Poincar\'e inequality holds for
$$\ll=\sup_{r>0} \Big\{K_0(r)\land \ff{K_1(r)}r\Big\}=\sup_{r>0} \Big\{\ff 1 r \land \Big(r_0-\ff m r\Big)\land \Big(Kr-r_0-\ff 4 r\Big)\Big\}.$$
Since $\ff 1 r$ is decreasing in $r>0$ with range $(0,\infty)$ while $\Big(r_0-\ff m r\Big)\land \Big(Kr-r_0-\ff 4 r\Big)$ is increasing in $r>0$ with range $(-\infty, r_0)$, $\ll$ is reached by a unique number $r_1>0$ such that
$$\ff 1{r_1}= \Big(r_0-\ff m {r_1}\Big)\land \Big(Kr_1-r_0-\ff 4 {r_1}\Big).$$ Then the value of  $\ll$ can be fixed by considering the following two situations:
\beg{enumerate}\item[A.] If $r_0-\ff m {r_1}\le Kr_1-r_0-\ff 4 {r_1}$, we have $\ff 1 {r_1}= r_0-\ff m{r_1}$ so that
$r_1= \ff{m+1}{r_0} $ and hence, $\ll= \ff{r_0}{m+1}.$
\item[B.] If $Kr_1-r_0-\ff 4 {r_1} <r_0-\ff m{r_1}$, then $\ff 1 {r_1}= Kr_1-r_0-\ff 4{r_1}$ so that $r_1= \ff{r_0+\ss{r_0^2+20K}}{2K}$ and
$\ll= \ff{2K}{r_0+\ss{r_0^2+20K}}$.\end{enumerate}

To prove the validity of the log-Sobolev inequality, we observe that
$$c_b\le 1+ 2 \Big\{m\land \Big(\ff{r_0^2}{4K}+4\Big)\Big\}=:c_0.$$ Moreover, by the Meyer diameter theorem (see \cite{BL} and references within), (\ref{CD1}) with $K>0$ implies that the intrinsic distance induced by $\bar\GG$ is bounded by a constant $D>0$. Noting that
$$\GG_b(f)(x,y)=tf_x^2(x,y)+t^2\bar\GG((x,\cdot))(y),$$ we conclude that
$$\rr_b^2((x,y),(x',y')) \le \ff{|x-x'|^2}t + \ff {D^2}{t^2},\ \ t>0, (x,y), (x',y')\in M.$$ Thus, for any $\ll>\ff {c_0}4 \big(\ge \ff{c_b}4\big)$,
$\mu(\e^{\ll \rr_b^2(o,\cdot)^2})<\infty$ holds for $o\in M$ and large $t>0$. Therefore, as indicated in the end of Subsection 3.2, the log-Sobolev inequality is valid.  \end{proof}

\subsection{Example B}  Consider the Gruschin operator $L= f_{xx} + x^{2l} f_{yy}$ on $M:=\R^2,$
  where $l\in\mathbb N$. We have  $$\GG^{(0)}(f,g)(x,y):= \GG(f,g)(x,y)= (f_xg_x)(x,y) + x^{2l}(f_yg_y)(x,y)$$ and
 $L= X^2+Y^2$ for   $X=\ff{\pp}{\pp x}, Y=x^l\ff{\pp}{\pp y}.$
When $l\ge 2$,
$\{X,Y, \nn_XY=lx^{l-1}\ff{\pp}{\pp y}, \nn_YX=0\}$ does not span the whole space for $x=0$. So, as explained in point (a) in the Introduction, (\ref{GCD'}) is invalid.
Let  $$\GG^{(i)}(f,g)(x,y)=x^{2(l-i)}(f_yg_y)(x,y),\ \ 1\le i\le l.$$ It is easy to see that $W(x,y):= 1+ x^2 +\ff{y^2}{1+x^{2l}}$ is a smooth compact function such that
\beq\label{W2}LW \le CW,\ \ \tt\GG(W)\le CW^2\end{equation}  holds for some constant $C>0$.

\beg{prp}\label{P4.2} There exist two constants $\aa,\bb>0$ depending only on $l$ such that $(\ref{GC})$ holds for
$$K_0(r_1,\cdots, r_l)= -\aa\sum_{i=1}^l \ff{r_{i-1}^{i-1}}{r_i^{i}},\ \ K_i(r_1,\cdots, r_l)= \bb r_{i-1}, \ \ 1\le i\le l, r_0=1,r_i>0.$$
Consequently, Proposition \ref{P3.1} holds for $b_i(0)= c_it^{2l-1+i},$ where
$$c_0=1, \ c_i= \ff{2\bb c_{i-1}}{2l-1+i},\ \ \ 1\le i\le l,$$ and
$$c_b= \sup_{r\in (0,t)} \Big\{(2l-1)r^{2(l-1)} + 2\aa \sum_{i=1}^l\ff{c_{i-1}^{i-1}}{c_i^i} r^{2(l-i)}\Big\}\le C_0 (1+t^{l-1}),\ \ t>0$$ for some constant $C_0>0.$ \end{prp}

\beg{proof} According to  Proposition \ref{P3.1} for $b_i(s)= c_i(t-s)^{2l-1+i},\ s\in [0,t], 0\le i\le l,$ it suffices to verify (\ref{GC}) for the desired $\{K_i\}_{0\le i\le l}$, which satisfy $$b_i'(s)+2 b_0(s) K_i\Big(\ff{b_1}{b_0},\cdots, \ff{b_l}{b_0}\Big)(s)=0,\ \ \ 1\le i\le l, s\in [0,t]$$ and
$$b_0'(s)+2 b_0(s) K_0\Big(\ff{b_1}{b_0},\cdots, \ff{b_l}{b_0}\Big)(s)=(2 l-1)(t-s)^{2(l-1)} + 2\aa \sum_{i=1}^l\ff{c_{i-1}^{i-1}}{c_i^i} (t-s)^{2(l-i)},\ \ s\in[0,t].$$
It is easy to see that at point $(x,y)\in \R^2$ and for $1\le i\le l$,
\beg{equation*}\beg{split} \GG_2(f)&= f_{xx}^2 + l(2l-1)x^{2(l-1)}f_y^2+ x^{4l}f_{yy}^2 + 2x^{2l} f_{xy}^2 + 4 l x^{2l-1}f_yf_{xy}- 2lx^{2l-1}f_xf_{yy},\\
\GG_2^{(i)}(f)&= (l-i) (2l-2i-1)x^{2(l-i-1)}f_y^2 + 4 (l-i) x^{2l-2i-1}f_yf_{xy} +x^{2(l-i)}f_{xy}^2 +x^{2(2l-i)}f_{yy}^2.\end{split}\end{equation*} So, for
$r_0=1$ and $r_i>0, 1\le i\le l,$
\beg{equation*}\beg{split} &\GG_2(f) +\sum_{i=1}^l r_i \GG_2^{(i)}(f) \\
&\ge f_y^2\sum_{i=0}^l r_i (l-i) (2l-2i-1)x^{2(l-i-1)} +f_{yy}^2\sum_{i=0}^l r_i x^{2(2l-i)} - 2l x^{2l-1}f_xf_{yy}\\
&\qquad  + f_{xy}^2 \Big\{2x^{2l} +\sum_{i=1}^l r_i x^{2(l-i)}\Big\} +f_yf_{xy} \sum_{i=0}^l 4 r_i (l-i) x^{2l-2i-1}\\
&\ge f_y^2\sum_{i=1}^{l} r_{i-1} (l+1-i) (2l-2i+1)x^{2(l-i)} - \ff{l^2}{r_1} f_x^2 -4 \sum_{i=0}^{l-1} \ff{r_i^2(l-i)^2x^{4l-4i-2}}{r_{i+1}x^{2(l-i-1)}}f_y^2\\
&\ge \ff { f_y^2} 2 \sum_{i=1}^l r_{i-1} (l+1-i) (2l-2i+1)x^{2(l-i)}  -\ff{l^2}{r_1}f_x^2\\
&\qquad - \ff 1 2 \sum_{i=1}^l \Big\{\ff{8r_{i-1}^2(l+1-i)^2 x^{2(l+1-i)}}{r_i} -r_{i-1} (l+1-i) (2l-2i+1)
x^{2(l-i)}\Big\} f_y^2 \\
&\ge \ff 1 2 \sum_{i=1}^l r_{i-1}(l+1-i)(2l-2i+1)x^{2(l-i)} \GG^{(i)}(f) -\ff{l^2}{r_1}f_x^2-\sum_{i=1}^l \ff{\aa_i r_{i-1}^{i-1}}{r_i^{i}}x^{2l} f_y^2\end{split}\end{equation*} holds for some constants $\aa_i>0, 1\le i\le l,$ where the last step is due to the fact that for constants $A_i,B_i>0$,
$$A_i x^{-2(i-1)}-B_i x^{-2i}\le \sup_{s>0} \big\{A_i s^{i-1} -B_i s^i\big\}=\ff{(i-1)^{i-1}A_i^i}{i^i B_i^{i-1}}.$$
\end{proof}

%Since when $l\ge 2$ (\ref{EX}) does not hold for the above $\{\GG^{(i)}\}_{1\le i\le l}$, results in Proposition \ref{P3.12} remain open for this example.

\subsection{Example C} Consider $Lf=f_{xx}+x^2f_{yy}+y^2f_{zz}=X_1^2+X_2^2+X_3^2$ on $\R^3$, where $X_1=\ff{\pp}{\pp x}, X_2= x\ff{\pp}{\pp y},
X_3= y\ff{\pp}{\pp z}$. It is easy to see that $\{X_i,\nn_{X_i}X_j\}_{1\le i\le 3}$ does not span $\R^3$ when $x=0$, so that (\ref{GCD'}) is not available according to observation (a) in Introduction. We have
$$\GG^{(0)}(f,g)(x,y,z):=\GG(f,g)(x,y,z)= (f_{x}g_x)(x,y,z) + x^2 (f_{y}g_y)(x,y,z)+y^2(f_zg_z)(x,y,z).$$
Let $$\GG^{(1)}(f,g)(x,y,z)=(f_yg_y)(x,y,z)+ x^2 (f_zg_z)(x,y,z),\ \ \GG^{(2)}(f,g)(x,y,z)= (f_zg_z)(x,y,z).$$ It is easy to see that
$W(x,y,z):=1+x^2+y^2$ is a smooth compact function on $\R^3$ such that (\ref{W2}) holds for some constant $C>0.$

\beg{prp}\label{P4.3}  $(\ref{GC})$ holds for
$$K_0(r_1,r_2)= -\Big(\ff 5 {r_1}+\ff{2 r_1}{r_2}\Big),\ \ K_1(r_1, r_2)=1-\ff{4r_1^2}{r_2}, \ \ K_2(r_1,r_2)=r_1,\ \  r_1,r_2>0.$$
Consequently, Proposition \ref{P3.1} holds for $$b_0(0)= t, \ \ b_1(0)= \ff {t^2}7,\ \ b_2(0)= \ff {2t^3}{21},$$ and
$c_b= 77.$   \end{prp}

\beg{proof} We first prove (\ref{GC}) for the desired $K_i, 0\le i\le 2.$ It is easy to see that at point $(x,y,z),$
\beg{equation*}\beg{split}& \GG_2(f)=f_{xx}^2 +f_y^2+x^2f_z^2+2x^2f_{xy}^2+2y^2f_{xz}^2+x^4f_{yy}^2+2x^2y^2f_{yz}^2+y^4f_{zz}^2\\
&\qquad\qquad+4xf_yf_{xy}+4x^2yf_zf_{yz}-2xf_xf_{yy} -2x^2yf_yf_{zz},\\
&\GG_2^{(1)}(f)=f_z^2+ f_{xy}^2+x^2f_{yy}^2+(y^2+x^4)f_{yz}^2+x^2f_{xz}^2+ x^2y^2f_{zz}^2+4xf_zf_{xz}-2yf_yf_{zz},\\
&\GG_2^{(2)}(f) = f_{xz}^2+x^2f_{yz}^2 +y^2f_{zz}^2.\end{split}\end{equation*}Therefore,
\beg{equation*}\beg{split} &\GG_2(f)+r_1\GG_2^{(1)}(f)+r_2\GG^{(2)}(f)\\
&\ge \GG^{(1)}(f)+ r_1\GG^{(2)}(f) +\big\{(2x^2+r_1)f_{xy}^2+4x f_yf_{xy}\big\}+\big\{(2y^2+r_1x^2+r_2)f_{xz}^2+4r_1xf_zf_{xz}\big\}\\
&\quad +\big\{(x^4+r_1x^2)f_{yy}^2-2xf_xf_{yy}\big\}+\big\{(2x^2y^2+r_1y^2+r_1x^4+r_2x^2)f_{yz}^2 +4 x^2yf_zf_{yz}\big\}\\
&\quad +\big\{(y^4+r_1x^2y^2+r_2y^2)f_{zz}^2 -2x^2yf_yf_{zz}-2r_1yf_yf_{zz}\big\}\\
&\ge \GG^{(1)}(f)+r_1\GG^{(2)}(f) -\ff{4x^2}{2x^2+r_1}f_y^2 -\ff{4r_1^2x^2}{r_2+r_1x^2+2y^2}f_z^2 -\ff{x^2}{x^4+r_1x^2}f_x^2\\
&\quad -\ff{4x^4y^2}{2x^2y^2 +r_1y^2+r_1x^4+r_2x^2}f_z^2 -\ff{(r_1+x^2)^2}{y^2+r_1x^2 +r_2}f_y^2\\
&\ge \GG^{(1)}(f)+r_1\GG^{(2)}(f) -\ff{4x^2}{r_1}f_y^2 -\ff{4r_1^2x^2}{r_2}f_z^2 -\ff{1}{r_1}f_x^2
 -\ff{4y^2}{r_1}f_z^2 -\ff{x^2}{r_1}f_y^2 -\ff{2r_1x^2}{r_2}f_y^2-\ff{r_1^2}{r_2}f_y^2 \\
 &\ge \Big(1-\ff{4r_1^2}{r_2}\Big)\GG^{(1)}(f) + r_1\GG^{(2)}(f)-\Big(\ff 5{r_1} +\ff{2r_1}{r_2}\Big)\GG(f).\end{split}\end{equation*}
This implies (\ref{GC}) for the claimed $K_i, 0\le i\le2 .$

Next, take
$$b_0(s)= t-s,\ \ b_1(s) =\ff 1 7 (t-s)^2,\ \ b_2(s)= \ff 2 {21}(t-s)^3,\ \ s\in [0,t].$$
Then
\beg{equation*}\beg{split} &\Big\{b_1' + 2b_0 K_1\Big(\ff{b_1}{b_0},\ff{b_2}{b_0}\Big)\Big\}(s) =-\ff{2(t-s)}7 +2(t-s)\Big(1-\ff 6 7\Big)=0,\\
&\Big\{b_2' + 2b_0 K_2\Big(\ff{b_1}{b_0},\ff{b_2}{b_0}\Big)\Big\}(s) =-\ff{6(t-s)}{21} +\ff{2(t-s)^2} 7  =0,\\
&\Big\{b_0' + 2b_0 K_0\Big(\ff{b_1}{b_0},\ff{b_2}{b_0}\Big)\Big\}(s) =-1-2(t-s)\Big(\ff{5}{t-s}+\ff 3{t-s}\Big) =-77.\end{split}\end{equation*}
Therefore, the second assertion holds.\end{proof}
%Again,   (\ref{EX}) does not hold  so that   Proposition \ref{P3.12} does not apply to this example.

\section{An extension of Theorem \ref{T1.1}}

As mentioned in the end of Section 1, for some highly degenerate subelliptic operators (\ref{GC}) is only available for   non-positively definite differential forms
$\GG^{(i)}$ and for some $r_1,\cdots, r_l>0$. For instance, when $L=\ff{\pp^2}{\pp x^2}+x\ff{\pp }{\pp y}$. one has $\GG(f)= f_x^2$ and
$\GG_2(f)= f_{xx}^2-f_xf_y.$ So, to verify (\ref{GC}), it is natural  to take $\GG^{(1)}(f,g) =-\ff {f_xg_y+f_yg_x}2,$ which is however not positively definite.
See Example 5.1 below for details.

To investigate such operators, we make use of the following   weaker version of assumption {\bf (A)}. We call bilinear symmetric form $\bar\GG:
C^3(M)\times C^3(M)\to C^2(M)$ a $C^2$ symmetric differential form, if
$$\bar\GG(fg,h)= f\bar\GG(g,h)+g\bar\GG(f,h),\ \ \ \bar\GG(f,\phi\circ g)=(\phi'\circ g)\bar\GG(f,g),\ \  \ f,g\in C^3(M), \phi\in C^1(\R)$$ holds.

\paragraph{(B)} \emph{There exist some $C^2$ symmetric differential form $\{\GG^{(i)}\}_{1\le i\le l},$ a non-empty set $\OO\subset (0,\infty)^l$, a smooth compact function $W\ge 1,$ and some  function $\{K_i\}_{0\le i\le l} \subset C(\OO)$ such that \beg{enumerate}\item[$(B1)$] $\GG_2+\sum_{i=1}^l r_i \GG_2^{(i)}\ge \sum_{i=0}^l K_i(r_1,\cdots, r_l)\GG^{(i)},\ \ (r_1,\cdots, r_l)\in\OO,$ where $\GG^{(0)}=\GG$.
\item[$(B2)$] $LW\le CW$ and $\sum_{i=0}^l |\GG^{(i)}(W)|\le CW^2$ hold for some constant $C>0.$
\item[$(B3)$] There exist $\vv>0$ and  $\tt r=(\tt r_1,\cdots,\tt r_l)\in\OO$ such that
$$\tt\GG:=\GG+\sum_{i=1}^l\tt r_i \GG^{(i)} \ge \vv \sum_{i=0}^l |\GG^{(i)}|.$$\end{enumerate} }

\beg{thm}\label{T5.1} Assume {\bf (B)}. For fixed $t>0$, let $\{b_i\}_{0\le i\le l}\subset C^1([0,t])$ be strictly positive on $(0,t)$ such that
\beg{enumerate}\item[$(i)$] $(\ff{b_1}{b_0},\cdots,\ff{b_l}{b_0})(s)\in \OO$ holds for all $s\in (0,t);$\item[$(ii)$]
$ b_i'(s) + 2 \big\{b_0 K_i (\ff{b_1}{b_0},\cdots,\ff{b_l}{b_0} )\big\}(s) = 0,\ \ s\in (0,t), 1\le i\le l.$\end{enumerate} Then assertions in $(1)$ and $(2)$ of Theorem $\ref{T1.1}$ hold. \end{thm}

\beg{proof} By $(B1)$ and $(B3)$, $\tt\GG_2\ge K\tt\GG$ and $\tt\GG\ge \vv \sum_{i=0}^l |\GG^{(i)}|$ hold for some $K\in\R$ and $\vv>0.$ Combining these with $(B2)$ and repeating the proof of Lemma \ref{L2.1}, we conclude that  $\{\GG^{(i)}(P_\cdot f)\}_{0\le i\le l}$ are bounded on $[0,t]\times M.$
Therefore, due to $(i)$ and $(ii)$ the proof of Theorem \ref{T1.1} works also for the present case. \end{proof}

To illustrate this result, we consider the following example mentioned in the beginning of this section, where the resulting gradient and Hanrack inequalities have the same   time behaviors as the corresponding ones presented in \cite[Corollaries 3.2 and 4.2]{GW} by using coupling methods.
In this example, it is easy to find correct choices of $W$, $\GG^{(i)}, K_i$ and $\OO$  such that assumption {\bf (B)} and condition (i) in Theorem \ref{T5.1} hold. The technical (also difficult)  point is to construct   functions $\{b_i\}_{i=0}^l$ such that condition (ii) holds and
$\sum_{i=0}^l b_i(0)\GG^{(i)}$ is an elliptic square field.

\paragraph{Example 5.1.} Consider $L=\ff{\pp^2}{\pp x^2}+x\ff{\pp }{\pp y}$ on $\R^2$. We have
$$\GG^{(0)}(f,g):=\GG(f,g)= f_xg_s.$$ Let
$$\GG^{(1)}(f,g)= -\ff 1 2(f_xg_y+f_yg_x),\ \ \ \GG^{(2)}=f_xg_y.$$ Then {\bf (B)} holds for $W(x,y)=1+x^2+y^2, \OO=\{(r_1,r_2): r_1,r_2>0, r_1^2\le 4r_2\},$ and
$$K_0(r_1,r_2)= 0, \ \ K_1(r_1,r_2)= 1+r_2,\ \ K_2(r_1,r_2)= \ff {r_1}2. $$ Moreover, $(\ref{EX})$ holds. Consequently, for any $\theta\in (\ff 3 2,2)$, letting
$t_\theta>0$ be the unique solution
$$(\coth\big[\ss 2 t_\theta\big] -1)^2 = 2\Big(\ss 2 \sinh \big[\ss 2 t_{\theta}\big] -2 t_\theta^2\Big),$$
there holds
\beq\label{LL5.1} \ff{2-\theta} 2 \bigg\{(t\land t_\theta) (P_tf)_x^2 +\ff {(t\land t_\theta)^3} 3 (P_tf)_y^2\bigg\} \le (P_t f)\big\{P_t(f\log f)- (P_tf)\log P_t f\big\},\ \ t>0\end{equation} so that the Harnack inequality
\beq\label{LLH}(P_t f)^\aa((x,y))\le (P_tf^\aa(x,y))\exp\bigg[ \ff {\aa}{2(2-\theta)(\aa-1)} \Big(\ff{|x-x'|^2}{t\land t_\theta}+\ff{3|y-y'|^2}{(t\land t_\theta)^3}\Big)\bigg]\end{equation}
holds for all $\aa>1,  t>0, (x,y),(x',y')\in\R^2$ and positive $f\in \B_b(\R^2)$.

\beg{proof} Obviously, $(B2)$ holds for the given $W$ and $(B3)$ holds for $r_1=r_2=1$ (hence $(r_1,r_2)\in\OO)$  and $\vv =\ff 1 4.$ Next, it is easy to see that (\ref{EX}) holds and
$$\GG_2(f)= f_{xx}^2-f_x f_y,\ \ \GG_2^{(1)}(f)= \ff 12 f_y^2 -f_{xx}f_{xy},\ \ \GG^{(2)}(f)= f_{xy}^2-f_xf_y.$$
Then, for $r_1,r_2>0$ with $r_1^2\le 4r_2$, we have
$$\GG_2(f)+r_1 \GG_2^{(1)}(f)+r_2\GG_2^{(2)}(f)\ge (1+r_2)\GG^{(1)}(f) +\ff{r_1}2 \GG^{(2)}.$$ Therefore, $(B1)$ holds.

To prove (\ref{LL5.1}) and (\ref{LLH}),  we take  $l=2$ and
$$b_0(s)= t-s,\ \ b_1(s)= \coth\big[\ss 2 (t-s)\big] -1,  \ \ b_3(s)=\ff 1 2 \sinh\big[\ss 2 (t-s)\big] -(t-s),\ \  s\in [0,t].$$ Then it is easy to see that $(ii)$ holds.    Observing that the function
$$\psi(r):= \ff{(\coth r-1)^2}{r\sinh r -r^2},\ \ r>0$$ is increasing with $ \lim_{r\to\infty}f(r)=\infty$ and $\lim_{r\to 0} f(r)=\ff 3 2,$ for any $\theta\in (\ff 3 2,2)$ the claimed quantity $t_\theta$ exists uniquely and for any $s\le t_\theta$,
\beq\label{LL*}\ff{b_1^2(s)}{b_0^2(s)}\le \ff{\theta^2 b_2(s)}{b_0(s)}\le \ff{4b_2(s)}{b_0(s)}.\end{equation} Thus, $(\ff{b_1}{b_0}, \ff{b_2}{b_0})\in\OO$ holds on $[0,t]$ provided  $t\le t_\theta.$
Since the right-hand side of (\ref{LL5.1}) is increasing in $t$, we may and do assume that $t\in (0, t_\theta].$
Noting that $b_0(0)=1$ and $b_i(t)-0$ for $0\le i\le 2$, it follows from Theorem \ref{T5.1} and (\ref{LL*}) that
\beg{equation*}\beg{split} &(P_t f) \big\{P_t(f\log f)- (P_t f)\log P_t f\big\}
 \ge\sum_{i=0}^2 b_i(0) \GG^{(i)} (P_tf) \\
 &=\ff{ (2-\theta)b_0(0)} 2 \Big\{(P_tf)_x^2+\ff{b_2(0)}{b_0(0)} (P_tf)_y^2\Big\}+ b_0(0)\Big\{\ff\theta 2 (P_tf)_x^2 -\ff{b_1(0)}{b_0(0)} (P_tf)_x(P_tf)_y
 +  \ff{\theta b_2(0)}{2b_0(0)}(P_tf)_y^2\Big\}\\
 &\ge \ff{2-\theta} 2 \Big\{t (P_tf)_x^2 +\Big(\ff{\sinh[\ss 2 t]}{\ss 2} -t\Big) (P_tf)_y^2\Big\}\ge \ff{2-\theta} 2 \Big\{t (P_tf)_x^2 +\ff{t^3} 3 (P_tf)_y^2\Big\}.\end{split}\end{equation*} Therefore, (\ref{LL5.1}) holds.

 Finally, letting $\bar\rr$ be   the intrinsic distance induced by the square field
$$ \bar\GG(f):= \ff{2-\theta} 2 \Big\{ (t\land t_\theta) (P_tf)_x^2 +\ff {(t\land t_\theta)^3}3 (P_tf)_y^2\Big\},$$ we have
$$\bar\rr((x,y),(x'y'))^2= \ff 2{2-\theta } \Big(\ff{|x-x'|^2}{t\land t_\theta}+\ff{3|y-y'|^2}{(t\land t_\theta)^3}\Big).$$
Then the desired Harnack inequality follows from Lemma \ref{L3.4} for $\gg(\dd)=\ff 1 {4\dd}$ since (\ref{LL5.1}) is equivalent to
$$\ss{\bar\GG(P_tf)}\le  \dd \big\{P_t(f\log f)- (P_tf)\log P_t f\big\} +\ff {P_tf} {4\dd},\ \ \dd>0.$$
\end{proof}

\end{document}